\documentclass[11pt]{amsart}
\title{On tolerances representable as $R \circ R^-$}
\author{Paolo Lipparini}
\address{Dipartimento di Matematica, Viale della Ricerca Scientifica,
II Universit\`a di Roma (Tor Verguccia),
 ROME 
ITALY}
\thanks{The author has received support from MPI and GNSAGA.
We acknowledge useful correspondence with G. Cz\'edli.} 

\email{lipparin@axp.mat.uniroma2.it}
\urladdr{http://www.mat.uniroma2.it/\textasciitilde lipparin}
\newtheorem{Theorem}{Theorem}
\newtheorem{proposition}[Theorem]{Proposition}

\newtheorem{theorem}[Theorem]{Theorem}

\newtheorem{corollary}[Theorem]{Corollary}

\theoremstyle{definition}

\newtheorem{definition}[Theorem]{Definition}

\newcommand{\alg}{\mathbf} 
\def\v{\mathcal V}  

\DeclareMathOperator{\Thetax}{\Theta} 
 
\DeclareMathOperator{\RR}{\mathit{R}}

\begin{document}

\begin{abstract} 
We give examples and counterexamples concerning varieties
in which every tolerance is representable as $R \circ R^-$,
for some reflexive and admissible relation $R$.
\end{abstract} 

\maketitle

In \cite{contol} we introduced the following definitions.

\begin{definition}\label{rep}   
A tolerance
$\Theta$ of some algebra ${\alg A} $  is {\em representable} 
if and only if there exists a compatible and reflexive
relation $R$ on ${\alg A} $ such that 
$\Theta= R \circ R^-$ (here, $R^-$ denotes the converse of $R$).

A tolerance
$\Theta$ of some algebra ${\alg A} $  is {\em weakly representable} 
if and only if there exists a set $K$ (possibly infinite) and there are
compatible and reflexive
relations $R_k$ ($k \in K$) on ${\alg A} $ such that 
$\Theta= \bigcap _{k \in K} (R_k \circ R_k^-)$. 
\end{definition}

The definitions are motivated by the following Theorem from \cite{contol}.

\begin{theorem}\label{contolthm} 
For every variety $\v$ and 
for every pair of terms $p$, $q$ (of the same arity) for the language
$ \{\circ, \cap\} $, if $p$ is regular, then 
the following are equivalent:

(i) $\v $ satisfies the congruence identity
$p(\alpha_1, \dots, \alpha_n) \subseteq q(\alpha_1, \dots, \alpha_n )$.

(ii) The tolerance identity
$p(\Theta_1, \dots, \Theta_n) \subseteq q(\Theta_1, \dots, \Theta_n)$
holds for every algebra $ {\alg A} $ in 
$\v$ and
for all representable tolerances
$\Theta_1, \dots, \Theta_n $ of $ {\alg A} $.

(iii) The tolerance identity
$p(\Theta_1, \dots, \Theta_n) \subseteq q(\Theta_1, \dots, \Theta_n)$
holds for every algebra $ {\alg A} $ in 
$\v$ and
for all weakly representable tolerances
$\Theta_1, \dots, \Theta_n $ of $ {\alg A} $.
 
(iv) $\v$ satisfies the tolerance identity
$p(\Theta_1 \circ \Theta_1, \dots, \Theta_n \circ \Theta_n) \subseteq 
q(\Theta_1 \circ \Theta_1, \dots, \Theta_n \circ \Theta_n)$.
\end{theorem}

We say that a term $p$ is regular
if and only if in the labeled graph associated with $p$
no vertex is adjacent with two distinct edges labeled
with the same name (see \cite{C2, czmsd, CD, contol} for details).

The aim of the present paper is to study the notion of a (weakly)
representable tolerance in more detail.

We first show that all tolerances in algebras without operations 
are weakly representable.

\begin{proposition} \label{sets}
If $ {\alg A} $ is an algebra belonging to the variety of sets
(that is, an algebra without operations) then every tolerance of
$ {\alg A} $ is weakly representable.
\end{proposition} 

\begin{proof} 
Let $ {\alg A} $ be an algebra without operations.
For every pair of distinct elements $a,b \in A$ let $ \Theta _{ab}$
be the reflexive and symmetric relation such that 
$(x,y) \in \Theta $ if and only if $ \{x,y\} \not= \{a,b\} $.

$ \Theta _{ab}$ is representable: define
$R$ by $x \RR y$ if and only if either $x=y=a$, or $x=y=b$,
or $x \not \in \{a,b\} $. $R$ is clearly reflexive, and is compatible since
$ {\alg A} $ has no operations. It is easy to see that
$ \Theta _{ab}=R \circ R^-$.

If $\Theta$ is any tolerance of $ {\alg A} $ then $\Theta$
is weakly representable, since $ \Theta= \bigcap _{(a,b) \not\in \Theta} \Theta _{ab}$.
 \end{proof}

In contrast to Proposition \ref{sets}, in algebras without
operations there can be non representable tolerances.
Such tolerances remain non representable if we add a certain kind of operations.

\begin{proposition} \label{sets2}
(i) In the 5-element algebra without operations there is a non representable tolerance.

(ii) There exists a 7-element semilattice with a non representable tolerance.

(iii) There exists a 7-element algebra with
a majority operation with a non representable tolerance
(a \emph{majority operation} is a ternary operation $f$ satisfying
$x=f(x,x,y)=f(x,y,x)=f(y,x,x)$).
\end{proposition} 

\begin{proof}
(i) Let $a,b_1,b_2,b_3,c$ denote the elements of the 5-element algebra
without operations, and let
$ \Theta$ be the smallest reflexive and symmetric relation
such that $ a \Thetax b_i $ and $b_i \Thetax c$ for $i=1,2,3$.

$ \Theta $ is a tolerance, since the algebra has no operations,
and it is easy to see that $ \Theta $ is not representable. Indeed, 
if $R$ is reflexive and $ \Theta = R \circ R^-$ then $ R \subseteq \Theta $
and $ R^- \subseteq \Theta $,
hence either $a \RR b_1$ or $b_1 \RR a$.   
Suppose that $a \RR b_1$ (the case $b_1 \RR a$ is similar).
If $c \RR b_1$ then $a \RR \circ \RR^- c$, that is, $a \Thetax c$,
which is false, hence necessarily $b_1 \RR c$. Continuing in the same way
we obtain both $b_2 \RR a$ and $b_3 \RR a$, which implies
$b_2 \RR \circ \RR^- b_3$, hence $b_2 \Thetax b_3$, contradiction.

(ii) Consider the semilattice $S$ with $6$ minimal elements
$a,b_1,b_2,b_3,b_4, c$ and with a largest element $1$.
Let $ \Theta $ be the smallest 
reflexive and symmetric relation
such that $1$ is $ \Theta $-related to all elements of $S$, and
such that
$ a \Thetax b_i $ and $b_i \Thetax c$ for $i=1,2,3,4$.

It is easy to check that $ \Theta $ is a tolerance.
Suppose by contradiction that $ \Theta $ is representable
as $ R \circ R^-$. If $x, y$  are minimal elements of $S$ and
both $x \RR 1$ and  $y \RR 1$, then $x \RR \circ \RR^- y$,
hence $x \Thetax y$.   
Thus 
$ | \{  x \in S | x \text{ is minimal and } x \RR 1\}  | \leq 2$,
since in $S$ there do not exist 3 pairwise $ \Theta $-connected minimal elements.

We can now repeat the arguments in (i) restricting
ourselves to minimal elements $x$ such that not $x \RR 1$.
 
(iii) Consider the lattice $ \langle L, +, \cdot \rangle$ with $6$ atoms
$a,b_1,b_2,b_3,b_4,c$ and with a largest element $1$
and a smallest element $0$.
If $f$ is the ternary operation defined by
$f(x,y,z)=(x+y)(x+z)(y+z)$ then
$\langle L\setminus \{0\}, f \rangle $ is an algebra, since
$L\setminus \{0\}$
 is closed under $f$.
We have that $f$ is a majority operation, 
and the same tolerance as in (ii)
is not representable.
\end{proof}

Even if we have showed that a majority term does not necessarily imply representability of
tolerances, we can show that lattices have representable tolerances.

\begin{proposition}\label{semilattices}
Suppose that the algebra 
${\alg A}$
has binary terms $\vee$ and 
$\wedge$ such that 
$\vee$ defines a join-semilattice operation,
the identities
$a \wedge (a \vee b)= a$,
$(a \vee b) \wedge b= b$
 are satisfied for all elements $a,b \in A$,
and the semilattice order induced by $\vee$ 
is a compatible relation on $ {\alg A} $.
Then all tolerances of 
${\alg A}$
are representable.

In particular, all tolerances in a lattice are representable.
\end{proposition}

\begin{proof}
If $ \Theta $ is a tolerance of $ {\alg A} $,  let $R= \Thetax \cap \leq$. 
$R$ is compatible since both $ \Theta $ and $\leq$ are compatible. 

If $a \Thetax b$ then $a = a \vee a \Thetax a \vee b$,
and $a \leq a \vee b$, thus $a \RR a \vee b$.
Similarly, $b \RR a \vee b$, that is, $a \vee b \RR^- b$, thus
$\Theta \subseteq R \circ R^-$.

Conversely,
if $(a,b) \in  R \circ R^-$,
say $a \RR c  \RR^- b$,
then 
$a \leq c$, thus
$c=a \vee c$, hence
$a= a \wedge (a \vee c)= a \wedge c$;
similarly,
$c \wedge b=b$,
hence
$a=a \wedge c \Thetax c \wedge b= b$,
since both $ R \subseteq \Theta $
and $ R^- \subseteq \Theta $.
Thus
$a \Thetax b$.
We have proved $ R \circ R^-\subseteq \Theta $.
\end{proof} 

We now proceed to show that if $ {\alg A} $ has a tolerance $ \Theta $
which is not a congruence, then we can add operations to $ {\alg A} $
in such a way that, in the expanded algebra, $ \Theta $ is not
even weakly representable. As a consequence, a Mal'cev condition
$\mathcal M$ implies that every tolerance is representable if and only if 
$\mathcal M$ implies congruence permutability (Corollary \ref{mal}).

\begin{proposition} \label{add}
Let $ {\alg A} $ be any algebra, and let $ \Theta$ be a tolerance of $ {\alg A} $.
Then there is an expansion $ {\alg A} ^+$ of $ {\alg A} $ by unary operations
such that $ \Theta$ is a tolerance of 
$ {\alg A} ^+$, and any non trivial reflexive compatible relation of $ {\alg A} ^+$ contains $ \Theta$.
\end{proposition} 

\begin{proof} 
Let $ {\alg A} ^+$ be obtained from $ {\alg A} $ by adding,
for every $a,b\in A$ such that $a \Thetax b$, 
and for every function $f: A \to \{a,b\} $,
a new unary operation which represents the function.
Since $a \Thetax b$, $\Theta$ is a tolerance of $ {\alg A} ^+$.

If $R$ is a non trivial reflexive compatible relation
of $ {\alg A} ^+$, there exist $c \not= d \in A$ such that $c \RR d$.
However, for every $a \Thetax b$ there is a function
such that $f(c) =a$ and $f(d)=b$, thus
$a=f(c) \RR f(d)=b$, since $R$ is compatible.
This proves that $R \subseteq \Theta$.
\end{proof}

\begin{corollary} \label{expans}
If $ {\alg A} $ is an algebra and $ \Theta$ is a tolerance of $ {\alg A} $
which is not a congruence,
then there is an expansion $ {\alg A} ^+$ of $ {\alg A} $ by unary operations
such that $ \Theta$ is a tolerance of 
$ {\alg A} ^+$ and $ \Theta $  is not representable in $ {\alg A} ^+$.
Actually, $ \Theta$ is not even weakly representable in $ {\alg A} ^+$.
  \end{corollary}

  \begin{proof} 
Let $ {\alg A} ^+$ be an expansion of $ {\alg A} $ as given by Proposition 
\ref{add}. 
$ \Theta$ is a tolerance of $ {\alg A} ^+$ by Proposition \ref{add}; 
moreover,  $ \Theta$ 
is non trivial, since the trivial tolerance is a congruence.
Suppose by contradiction that
 $ \Theta= R \circ R^-$ for some reflexive and admissible
relation $R$ on $ {\alg A} ^+$, hence $R$ and $R^-$ are non trivial, thus
$R \supseteq \Theta$ and $R^- \supseteq \Theta$,
by Proposition \ref{add}. Then
$ \Theta= R \circ R^- \supseteq \Theta \circ \Theta$,
and this implies that $ \Theta$ is a congruence of $ {\alg A} ^+$,
hence a congruence of $ {\alg A}$, contradiction.
The proof that $\Theta$ is not weakly representable 
in $ {\alg A} ^+$ is similar.
\end{proof} 

The following result is probably known, but we give a proof, since
we know no reference for it.

\begin{proposition} \label{perm} 
(a) If $ {\alg A} $ is an algebra, and every tolerance of $ {\alg A} $
is a congruence, then all congruences of $ {\alg A} $ permute.

(b) A variety $\v$ is congruence permutable if and only if 
every tolerance of every algebra in $\v$ is a congruence. 
   \end{proposition}

   \begin{proof} 
(a) If $ \alpha , \beta $ are congruences of $ {\alg A} $,
let $ \overline{ \alpha \cup \beta }$ denote the smallest tolerance containing
$ \alpha $ and $ \beta $, which is the smallest admissible relation
containing $ \alpha  \cup \beta  $.
Notice that $ \overline{ \alpha \cup \beta } \subseteq \beta \circ \alpha $.

 By assumption, $\overline{ \alpha \cup \beta }$
is a congruence.
Then $ \alpha \circ \beta \subseteq \overline{ \alpha \cup \beta } 
\circ \overline{ \alpha \cup \beta }=\overline{ \alpha \cup \beta } \subseteq 
\beta \circ \alpha $.

(b) is immediate from (a) and the well known result that in 
permutable varieties every reflexive and admissible relation is a congruence
(see \cite{HM}, \cite[Proposition 143]{Sm}).
\end{proof}

Trivially, every congruence $ \alpha $ is representable, since
$ \alpha = \alpha \circ \alpha $. 
By Proposition \ref{perm}(b), 
congruence permutability, for varieties, implies that
every tolerance is representable.
The next result shows that if a Mal'cev
condition $\mathcal M$ implies that
every tolerance is representable, then 
$\mathcal M$ implies congruence permutability.

\begin{corollary} \label{mal}
Let 
$\mathcal M$
be either a Mal'cev condition,
or a weak Mal'cev condition, or a strong Mal'cev condition.
The following are equivalent:

(i) $\mathcal M$ implies congruence permutability.

(ii) $\mathcal M$ implies that
every tolerance is representable.

(iii) $\mathcal M$ implies that
every tolerance is weakly representable.
\end{corollary}

\begin{proof} 
(i) $\Rightarrow$ (ii). Suppose that (i) holds. If $\v$ satisfies
$\mathcal M$,
then, by Proposition \ref{perm}(b), every tolerance 
in every algebra in $\v$ is a congruence, hence is representable. 
Thus, (ii) holds.

(ii) $\Rightarrow$ (iii) is trivial.

We shall prove (iii) $\Rightarrow$ (i) by contradiction.

Suppose that (i) fails.
Then there exists some variety  $\v$ which satisfies
$\mathcal M$
but which is not congruence permutable.
 By Proposition \ref{perm}(b), there is an algebra $ {\alg A} \in \v$
with a tolerance $ \Theta $ which is not a congruence. By Corollary 
\ref{expans}, $ {\alg A} $ can be expanded to an algebra
$ {\alg A} ^+$ in which $ \Theta $ is a tolerance which is not weakly representable.
By well known properties of Mal'cev conditions, the variety
generated by $ {\alg A} ^+$ still satisfies $\mathcal M$, and this contradicts
(iii).
\end{proof}

\begin{corollary} \label{mal2}
(i) The class of varieties $\v$ such that every tolerance in every
algebra in $\v$ is representable cannot be characterized by a weak Mal'cev condition.

(ii) The class of varieties $\v$ such that every tolerance in every
algebra in $\v$ is weakly representable cannot be characterized by a weak Mal'cev condition.
  \end{corollary} 

\begin{proof} 
If any of those classes could be characterized by some weak Mal'cev condition
$\mathcal M$, 
then, by Corollary \ref{mal}, $\mathcal M$ would imply
permutability.
This is a contradiction, since 
Propositions \ref{sets} and \ref{semilattices} 
provide examples of non permutable varieties 
in which every tolerance is (weakly) representable.
\end{proof}

\def\cprime{$'$} \def\cprime{$'$}

\providecommand{\bysame}{\leavevmode\hbox to3em{\hrulefill}\thinspace}

\end{document}